\documentclass{amsart}

\usepackage{amssymb}
\usepackage[all]{xy}
\usepackage{hyperref}

\newtheorem{thm}{Theorem}

\newtheorem{lem}[thm]{Lemma}
\newtheorem{cor}[thm]{Corollary}

\newtheorem{prop}[thm]{Proposition}

\newtheorem{conj}[thm]{Conjecture}
   
\theoremstyle{definition}
\newtheorem{defn}[thm]{Definition}

\newtheorem{say}[thm]{}
\newtheorem{exmp}[thm]{Example}

\newtheorem{rem}[thm]{Remark}          

\newtheorem*{ack}{Acknowledgments}      

\newtheorem{defn-thm}[thm]{Definition--Theorem}  
\newtheorem{defn-lem}[thm]{Definition--Lemma}  

\theoremstyle{remark}

\let \cedilla =\c
\renewcommand{\c}[0]{{\mathbb C}}  

\renewcommand{\o}[0]{{\mathcal O}} 
\newcommand{\z}[0]{{\mathbb Z}}

\renewcommand{\r}[0]{{\mathbb R}}

\newcommand{\dd}[0]{{\mathbb D}}

\newcommand{\p}[0]{{\mathbb P}}

\newcommand{\q}[0]{{\mathbb Q}}
\newcommand{\map}[0]{\dasharrow}

\newcommand{\spec}[0]{\operatorname{Spec}}
\newcommand{\pic}[0]{\operatorname{Pic}}

\newcommand{\rank}[0]{\operatorname{rank}}

\newcommand{\discrep}[0]{\operatorname{discrep}}
\newcommand{\totaldiscrep}[0]{\operatorname{totaldiscrep}}

\newcommand{\supp}[0]{\operatorname{Supp}}    
\newcommand{\red}[0]{\operatorname{red}}    
\newcommand{\codim}[0]{\operatorname{codim}}    
    
\newcommand{\proj}[0]{\operatorname{Proj}}

\newcommand{\coker}[0]{\operatorname{coker}}

\newcommand{\ex}[0]{\operatorname{Ex}}    
\newcommand{\diff}[0]{\operatorname{Diff}}

\newcommand{\chow}[0]{\operatorname{Chow}}

\newcommand{\hilb}[0]{\operatorname{Hilb}}

\newcommand{\rdown}[1]{\lfloor{#1}\rfloor}

\newcommand{\tsum}[0]{\textstyle{\sum}}

\newcommand{\coeff}[0]{\operatorname{coeff}}

\def\into{\DOTSB\lhook\joinrel\to}

\def\loccoh#1.#2.#3.#4.{H^{#1}_{#2}(#3,#4)}

\DeclareMathAlphabet{\mathchanc}{OT1}{pzc}%
                                {m}{it}

\usepackage[all]{xy}\xyoption{dvips}

\begin{document}
\bibliographystyle{amsalpha}

  \title{Moishezon morphisms}
  \author{J\'anos Koll\'ar}

\begin{abstract} We try to understand which morphisms of complex analytic spaces come from algebraic geometry. We start with a series of conjectures, and then give some partial solutions.
\end{abstract}

 \maketitle

\tableofcontents

 A proper, irreducible, reduced  analytic space  $X$ is {\it Moishezon} if it is
bimeromorphic to a  projective variety  $X^{\rm p}$,  and a 
   proper morphism of analytic spaces  $g:X\to S$   is {\it Moishezon} if it  is bimeromorphic to a projective
morphism  $g^{\rm p}:X^{\rm p}\to S$; see (\ref{mois.def.1}) and (\ref{mois.morph.defn}--\ref{loc.mois.morph.defn}) for details.

The aim of this note is to discuss a series of   questions about Moishezon morphisms, and give partial solutions to some of them.

We start with a list of conjectures in Section~\ref{sec.1}.
 Sections~\ref{sec.2}--\ref{sec.3} are mostly review; new results are in 
Sections~\ref{sec.4}--\ref{sec.7}.

\section{Open questions}\label{sec.1}

The theory of Moishezon spaces can be viewed as a special chapter  of the theory of algebraic spaces  (and later stacks). 
However, a deformation of a Moishezon space need not be Moishezon, thus we get a theory that is not algebraic. 
The  question we consider is the following.

\medskip
\begin{itemize}
\item  Which morphisms of complex amalytic spaces come from algebraic geometry, up to bimeromorphism?
\end{itemize}
\medskip

The main open problem may be Conjecture~\ref{main.question} and its special case Conjecture~\ref{main.question.1d}.  I have very little evidence  supporting them and several experts thought that they are likely wrong.

Conjectures~\ref{lim.of.M.M.conj}--\ref{main.question.gt} about the deformation theory of Moishezon spaces go back at least to Hironaka's unpublished thesis.

Conjecture~\ref{proj.mod.conj.c/lt} can be viewed as a   geometric form of deformation invariance of plurigenera \cite[Thm.1.2]{rao-tsa}, see also (\ref{def.i.pg.say}).

\begin{conj} \label{main.question}  
A proper morphism of analytic spaces  $g:X\to S$ is locally  Moishezon (\ref{loc.mois.morph.defn}) iff every irreducible component of every fiber is  Moishezon.
\end{conj}

{\it Comments \ref{main.question}.1.} We  check in (\ref{moi.fib.moi.cor}) that the fibers of a  Moishezon morphism are  Moishezon.
If $g$ is smooth, a positive answer is  in \cite[Thm.1.4]{rao-tsa-1} and \cite[Thm.1.2]{rao-tsa}; see (\ref{RT.vb.Zd.M.cor}).

If $\pi_1(S)$ is finite and $S$ is either Stein or quasi-projective, then  maybe $g$ is  also globally Mosihezon.
Easy examples  (\ref{fib.bdl.not.m}) show that  finiteness of   $\pi_1(S)$ is necessary   for the global variant.

 Note that we do not assume that $g$ is flat or that $X,S$ are smooth.
It is, however, quite likely that the above generality does not matter.
The semi-stable reduction theorem \cite{Abramovich-Karu00} suggests that there is a
projective, bimeromorphic morphism $S'\to S$ such that the main component of
$X\times_SS'$ is  bimeromorphic  to a  morphism $g':X'\to S'$ such that
$g'$ is flat, toroidal and 
$S'$ is smooth.

A positive  answer  for $g':X'\to S'$
does not automatically imply a  
positive  answer  for $g:X\to S$, but the method may generalize.

It is  reasonable to start with the case when $g$ is flat with mildly singular fibers.   The following may well be the key special case, where $\dd$ denotes the complex disc.

\begin{conj} \label{main.question.1d} Let $X$ be a smooth analytic space and
 $g:X\to \dd$  a proper morphism. Assume that the fiber $X_s$ is 
 Moishezon for $s\neq 0$, and $X_0$ is a simple normal crossing divisor whose irreducible components are  Moishezon.
Then  $g$  is Moishezon. 
\end{conj}

For an arbitrary proper morphism, the set of Mosihezon fibers need not be closed (\ref{bad.fams.exmp}), but the following could be true.

\begin{conj} \label{lim.of.M.M.conj}
Let $g:X\to \dd$ be a flat, proper morphism. Assume that
$X_0$ is irreducible with rational singularities  and $X_s$ is Moishezon for $s\neq 0$. 
Then $X_0$ is  Moishezon.
\end{conj}

{\it Comments \ref{lim.of.M.M.conj}.1.} There are 3 preprints \cite{pop09, barlet2017gauduchons, pop19} 
claiming a positive answer if  $X_0$ is smooth. 
These lie outside my expertise, but my understanding is that  not everyone is able to follow the arguments in them.

The analogous question for surfaces with cusp singularities has a negative answer; see (\ref{inoue.exmp}). Cusps are the simplest non-rational surface singularities. This suggested that either log terminal or  rational singularities may be the right class here.

\begin{conj} \label{main.question.gt} Let $X$ be a smooth analytic space and
 $g:X\to \dd$  a proper morphism. Assume that  one of the irreducible components of $X_0$ is of general type.
Then $g$ is Moishezon  and all other fibers are of general type (over a possibly smaller disc).
\end{conj}

{\it Comments \ref{main.question.gt}.1.} 
Smooth, projective K3 and elliptic surfaces have deformations that are
not even Moishezon, so  general type may be the best one can hope for.
We can harmlessly assume that $X_0$ is a reduced, simple normal crossing divisor.

If $X_0$ is irreducible  and smooth, this is posed in \cite[p.201]{sundararaman}; which in turn builds on  problems and conjectures in
\cite{iitaka,  MR0425189, nak-par, ueno}.

If $X_0$ is irreducible,  projective and has canonical singularities,
a positive answer is given  in \cite{k-defgt}.
Note, however, that  smooth Moishezon spaces can  have  unexpected 
deformations; see \cite{MR1107661, leb-poo}.

\medskip

Let  $g:X\to S$   be a   Moishezon morphism.  By definition, it  is bimeromorphic to a projective
morphism  $g^{\rm p}:X^{\rm p}\to S$. Thus the fibers  $X_s$ and $X^{\rm p}_s$ are 
bimeromorphic to each other for general $s\in S$,  but may be quite different for special points $s\in S$. The following conjecture suggests that, over 1-dimensional bases,   one can arrange   $X_s$ and $X^{\rm p}_s$ to be  
bimeromorphic to each other for every $s\in S$.

\begin{conj}\label{proj.mod.conj.c/lt}
 Let $g:X\to \dd$ be a flat, proper,   Moishezon morphism. Assume that  $X_0$ has canonical (resp.\ log terminal)  singularities.
 Then  $g$ is  
 fiberwise birational (\ref{fib.w.bim.defn}) to a  flat,   projective morphism
  $g^{\rm p}: X^{\rm p}\to \dd$ such that
\begin{enumerate}
\item  $X_0^{\rm p}$ has  canonical (resp.\ log terminal)  singularities,
\item  $X_s^{\rm p}$ has   terminal  singularities for $s\neq 0$, and
\item $K_{X^{\rm p}}$ is $\q$-Cartier.
\end{enumerate}
\end{conj}

{\it Comments \ref{proj.mod.conj.c/lt}.4.} 
This is where singularities inevitably enter the picture.
Even if  $g$ is a smooth family of projective surfaces, $X^{\rm p}$ may need to be singular; see for example \cite[Exmp.4]{k-defgt}.

If $g$ is smooth and the fibers are of general type, then
\cite[Thm.1.2]{rao-tsa} implies that the canonical models of the fibers
give the optimal choice for $g^{\rm p}: X^{\rm p}\to \dd$. (Here all fibers can have canonical  singularities.)

We give  a positive answer to the log terminal case,
 provided
$X_0$ is not uniruled, see (\ref{proj.mod.thm.lt}). The canonical case is discussed in 
 (\ref{can.main.question.1p.say}).

\begin{rem} My aim is to understand how much of the theory of Moishezon spaces fits into algebraic geometry, and especially minimal model theory. The paper \cite{rao-tsa} gave the impetus to try to organize this into  a systematic series of questions. 

Campana pointed out that several of these questions have analogs for
compact spaces of Fujiki's class $C$, and have a positive answer if we assume that the total space is of class $C$; see \cite{MR620706}.

A very different direction studies the place of the Moishezon property in the theory of compact complex manifolds.
Solutions of Conjectures~\ref{main.question}--\ref{main.question.1d}
are more likely to come from this approach.
 See \cite{pop11}  for a survey.
\end{rem}

\begin{ack} I thank D.~Abramovich, F.~Campana, J.-P.~Demailly, T.~Murayama, V.~Tosatti and C.~Xu
for  helpful comments.    I am grateful to 
S.~Rao and I-H.~Tsai for detailed remarks,  corrections and references.
Partial  financial support    was provided  by  the NSF under grant number
DMS-1901855.
\end{ack}

\section{Moishezon spaces}\label{sec.2}

We give a quick review of the theory of Moishezon spaces.

\begin{defn} \label{mois.def.1}
A proper, irreducible, reduced  analytic space  $X$ is {\it Moishezon} if it is
bimeromorphic to a  projective variety  $X^{\rm p}$. That is, there is a closed, analytic subspace $\Gamma\subset X\times X^{\rm p}$ such that the coordinate projections  $\Gamma\to  X$ and $\Gamma\to  X^{\rm p}$  are isomorphisms
on Zariski  open dense sets.

By Chow's theorem, any 2 such $X^{\rm p}$ are birational to each other, so $X$ acquires a unique algebraic structure.

 A proper analytic space $X$ is  Moishezon iff the 
irreducible  components of $\red X$ are  Moishezon\footnote{This is  not standard terminology.}.
Thus $X$ is Moishezon iff every irreducible component of its normalization is Moishezon.

\end{defn}

\begin{say}[Basic theorems] \label{basic.thms.say}
Let $X$ be a proper Moishezon space. 

\begin{enumerate}
\item There is a projective variety $X'$ and a bimeromorphic morphism
$X'\to X$  (Chow  lemma).  
\item For every $x\in X$ there is a pointed quasi projective scheme
$(x', X')$ and an \'etale morphism  $(x', X')\to (x, X)$.
\item If $X$ is normal then there is 
 a proper variety $Y$ and  a finite group $G$ acting on $Y$ such that $X\cong Y/G$. (Note that usually $Y$ can not be chosen projective.)
\item If $Z\to X$ is finite, then $Z$ is Moishezon. 
\item If $X\to Y$ is  surjective, then $Y$ is Moishezon.
\item Assume that $X$ is smooth. Then the usual Hodge decomposition 
$H^i(X, \c)=\oplus_{p+q=i} H^p(X, \Omega_X^q)$ holds.
\item $\hilb(X)$ and $\chow(X)$ are algebraic spaces whose irreducible components are proper (but the connected components may have infinitely many irreducible components). The connected components of the space of divisors $\chow_{n-1}(X)$ are proper. 
\item If $X$ has rational singularities then it is projective iff it is K\"ahler.
\end{enumerate}

{\it Hints of proofs.} Note that (1) is not obvious \cite{Moi-66}.
It also follows from the more general results of \cite{MR0393556}, and
one can easily modify the arguments in \cite[Tag 088U]{stacks-project}; the key step is probably 
\cite[Tag 0815]{stacks-project}.

(2) is quite hard; see  \cite{artin}.

For (3), cover $X$ with finitely many   $X'_i\to X$ as in (2). Then normalize $X$ in the Galois closure of the field extensions
$\c(X'_i)/\c(X)$.  One can then use this to get $Z\to X$ as a quotient of a finite morphism $Z'\to X'$ to obtain (4).

For (5), using (4) and (1) we may assume that $X\to Y$ is generically finite, say of degree $d$, and $X$ is projective. Then $y\mapsto [g^{-1}(y)]\in S^dX$ gives a bimeromorphic embedding
of $Y$ into the $d$th symmetric power of $X$. 

By direct computation, the existence of a Hodge decomposition is invariant under smooth blow ups, thus we get (6). A better argument is in \cite[Prop.1.3]{ueno-83}. 

For (7) see \cite{Artin69b, MR0399503, MR583821, fuj-c} and \cite[Sec.I.5]{rc-book}.

The  smooth case of (8) is proved in \cite{Moi-66}, the singular one in 
\cite{nam-02}. 
\medskip

{\it Remark.} The complements of closed analytic subsets form the open subsets of the Zariski topology. Note, however, that 2 open subsets can be biholomorphic to each other even if they are not birational. This is the main reason why one usually does not define `Moishezon' for non-proper anaytic spaces.

\end{say}

\section{Moishezon morphisms}\label{sec.3}

\begin{defn}[Projective morphisms] \label{proj.morph.defn}  A  proper morphism of analytic spaces  $g:X\to S$ 
 is  {\it projective} if $X$ can be embedded into  ${\mathbf P}_S:=\p^N\times S\to S$ for some $N$. Note that   some authors allow ${\mathbf P}_S\to S$ to be any (locally trivial) $\p^N$-bundle.
The 2 versions are  equivalent if   $S$ is Stein  or quasi-projective (the cases we are mostly interested in) but not in general.
\end{defn}

\begin{defn}[Moishezon morphisms] \cite{MR0369746, fuj-c} \label{mois.morph.defn}  
Assume now that $S$ is reduced. A  proper morphism of analytic spaces  $g:X\to S$   is {\it Moishezon} iff the
following equivalent conditions hold. 
\begin{enumerate}
\item  $g:X\to S$ is bimeromorphic to a projective
morphism  $g^{\rm p}:X^{\rm p}\to S$. That is, there is a closed subspace
$Y\subset X\times_S X^{\rm p}$ such that the coordinate projections
$Y\to X$ and $Y\to X^{\rm p}$ are bimeromorphic.
\item  There is a projective morphism of algebraic varieties
$G:{\mathbf X}\to {\mathbf S}$ and a meromorphic
$\phi_S:S\map {\mathbf S}$ such that
$X$ is bimeromorphic to   ${\mathbf X}\times_{\mathbf S}S$.
\end{enumerate}
Here (2) $\Rightarrow$ (1) is clear. To see the converse, note that   $g^{\rm p}:X^{\rm p}\to S$ is flat over a dense, open subset $S^\circ\subset S$, thus we get a meromorphic map
$\phi: S\map\hilb(\p^N)$. The pull-back of the universal family over
$\hilb(\p^N)$ is then bimeromorphic to $X$.

{\it Comment.} This is the right notion if $S$ is Stein or quasi projective, but, as with projectivity, there are different versions in general.
\medskip

Assume that $X$ is normal and the maps
$$
X\stackrel{\phi}{\map} X^{\rm p} \stackrel{\iota}{\into} {\mathbf P}_S
\eqno{(\ref{mois.morph.defn}.3)}
$$
show that $X\to S$ is Moishezon. Then $\iota\circ \phi:X\map {\mathbf P}_S$
is defined outside a codimension 2 closed subset, and
$(\iota\circ \phi)^*\o_{{\mathbf P}_S}(1)$ extends to a 
rank 1 reflexive sheaf $L$ on $X$. This $L$  `certifies' that $X$ is  Moishezon. This gives another equivalent characterization
(in case $X$ is normal, and $S$ is Stein or quasi projective.)
\begin{enumerate}\setcounter{enumi}{3}
\item  There is a rank 1, reflexive sheaf $L$ on $X$ such that 
the natural map
$X\map \proj_S(g_*L)$
is  bimeromorphic onto the closure of its image. 
\end{enumerate}
We call such a sheaf $L$  {\it very big}  (over $S$)\footnote{{\it Very} big is  not standard terminology, but it matches {\it very} ample.}.   
 $L$  is {\it  big}  (over $S$) if $L^{[m]}$ is very big for fome $m>0$, where 
$L^{[m]}$ denotes the reflexive hull of the $m$th tensor power.

Note that $L$ is  big (resp.\ very big) on $X\to S$ iff it is  big (resp.\ very big) on $X^\circ\to S^\circ$ on some dense, Zariski open $S^\circ\subset S$. 

{\it Warning.} By contrast it can happen that $X^\circ\to S^\circ$ is
Moishezon but  $X\to S$ is not, since the $L^\circ$ that certifies
Moishezonness need not extend to $X$; see (\ref{inoue.exmp}). 
\end{defn}

\begin{defn}[Locally Moishezon morphisms] \cite{MR0369746}  \label{loc.mois.morph.defn}  
 A   proper morphism of analytic spaces  $g:X\to S$    {\it locally Moishezon} if $S$ is covered by (Euclidean) open sets $S_i\subset S$ such that  each $g^{-1}(S_i)\to S_i$ is 
Moishezon.

{\it Comment.} This follows  standard usage of `locally' in algebraic geometry and it  works best for the purposes of 
Conjecture~\ref{main.question}.  However, it is {\em not} equivalent to  
the definition in
\cite{fuj-c}.
\end{defn}

\begin{exmp} Let $g:X\to S$  be a  proper morphism of analytic spaces, $S$  
 Moishezon. Then $g$ is Moishezon iff $X$ is Moishezon.
\end{exmp}

\begin{exmp} \label{fib.bdl.not.m}
Let $Z$ be a normal, projective variety with discrete automorphism group. Let $g:X\to S$ be  a fiber bundle with fiber $Z$ over a connected base $S$.
Then $g$ is Moishezon  $\Leftrightarrow$ $g$ is projective  $\Leftrightarrow$ the monodromy is finite.

There are rational and K3 surfaces with  infinite, discrete automorphism group. These lead to fiber bundles over the punctured disc  $\dd^\circ$ 
that are locally Moishezon but not globally Moishezon.
\end{exmp}

\begin{exmp}\label{inoue.exmp}
\cite{MR632841} studies examples where $X_0$ is an Inoue surface (which is  not Moishezon) with a  cusp  (which is log canonical), yet $X_s$ is a smooth  rational surface for $s\neq 0$. 
\end{exmp}

  Next we look at fibers of Moishezon morphisms.

\begin{lem} \label{ex.moi.lem}
Let  $g:X\to S$  be a proper, generically finite, dominant morphism of   normal, complex, analytic spaces.
Then $\ex(g)\to S$ is Moishezon.
\end{lem}

Proof. We prove the special case when the  smooth locus of $S$ is dense in $g\bigl(\ex(g)\bigr)$. This is a harmless assumption if $S$ is Stein (or quasi-projective), since we can compose $g$ with a finite
$S\to \c^{\dim S}$ (or with a quasi-finite
$S\to \p^{\dim S}$). A more heavy handed approach, which works in general,  is to
use a resolution  $S'\to S$ and replace $X$ by  the normalization of the main component of $X\times_SS'$. 

 Let $E_0$ be a $g$-exceptional divisor.
Set $(g_0:X_0\to S_0):=(g:X\to S)$ and $Z_0:=g_0(E_0)$.

If $g_i:X_i\to S_i $ and $E_i\subset X_i$ are already defined, we set
$Z_i:=g_i(E_i)$. Let $S_{i+1}$ be the normalization of the blow-up
$B_{Z_i}S_i$, $g_{i+1}:X_{i=1}\to S_{i+1}$ the normalization of the graph of
$X_i\to S_i\map S_{i+1}$ and $E_{i+1}\subset X_{i+1}$ the bimeromorphic transform of $E_i$.  (Note that $X_{i+1}\to X_i$ is a isomorphism over an open subset of $E_i$.)

Let $a(E_i, S_i)$ denote the vanishing order of the Jacobian of $g_i$ along $E_i$.  By an elementary computation   we get that
$$
a(E_{i+1}, S_{i+1})\leq a(E_i, S_i)+1-\codim(Z_i\subset S_i).
$$
Thus eventually we reach the situation when $\codim(Z_i\subset S_i)=1$,
hence $E_i\to Z_i$ is generically finite. 

Note that each $Z_{i+1}\to Z_i$ is projective, thus $E_i\to Z_0$ is Moishiezon by (\ref{basic.thms.say}.4), and so is $E_0\to S$. \qed

 \medskip

The following is the easy direction of  Conjecture~\ref{main.question}.

\begin{cor}\label{moi.fib.moi.cor} The fibers of a proper, Moishezon morphism are Moishezon.
\end{cor}

Proof. Let $g:X\to S$ be a proper, Moishezon morphism. It is
bimeromorphic to a projective morphism   $X^{\rm p}\to S$. We may assume $X^{\rm p}$ to be normal. Let $Y$ be the normalization of the closure of the graph of
$X\map X^{\rm p}$. 

Fix now $s\in S$. Let $Z_s\subset X_s$ be an irreducible component and
$W_s\subset Y_s$ an irreducible component that dominates $Z_s$. By (\ref{basic.thms.say}.5)
it is enough to show that $W_s$ is Moishezon. 

If $\pi: Y\to X^{\rm p}$ is generically an isomorphism along $W_s$, then
$W_s$ is bimeromorphic to an irreducible component of $X^{\rm p}_s$, hence Moishezon. Otherwise    $W_s\subset \ex(\pi)$. 
Now $\ex(\pi)\to  X^{\rm p}$ is Moishezon by (\ref{ex.moi.lem}) and
$\dim \ex(\pi)<\dim Y=\dim X$. 
Now $W_s$ is contained in a fiber of $\ex(\pi)\to S$, hence Moishezon by
 induction on the dimension.\qed

\begin{rem}\label{mois.morph.fib}
More generally, if $g:X\to S$ is  proper and Moishezon
and $T\to S$ is a morphism of analytic spaces  
then  $X\times_ST\to T$ is also  proper and  Moishezon. 
\end{rem}

The rest of this section  is a study of the set of Moishezon fibers for arbitrary proper morphisms of analytic spaces.
It is mostly a summary of some of the results of \cite{rao-tsa}, with occasional  changes.

\begin{defn} Let $g:X\to S$ be a  proper morphism of normal analytic spaces and  $L$  a line bundle on $X$. Set
\begin{enumerate}
\item $\operatorname{VB}_S(L):=\{s\in S:  L_s \mbox{ is very big on } X_s\}\subset S$,
\item  $\operatorname{GT}_S(X):=\{s\in S:  X_s \mbox{ is of general type}\}\subset S$,
\item  $\operatorname{MO}_S(X):=\{s\in S:  X_s \mbox{ is Moishezon}\}\subset S$,
\item  $\operatorname{PR}_S(X):=\{s\in S:  X_s \mbox{ is projective}\}\subset S$.
\end{enumerate}
\end{defn}

\begin{lem} \label{L.vb.Zd.M}
Let $g:X\to S$ be a  proper morphism of normal, irreducible analytic spaces and  $L$  a line bundle on $X$. Then $\operatorname{VB}_S(L)\subset S$ is
\begin{enumerate}
\item either  nowhere dense (in the analytic Zariski topology),
\item or  it contains a dense open subset of $S$, and
$g:X\to S$ is Moishezon.
\end{enumerate}
\end{lem}

Proof.  By passing to an open subset of $S$, we may assume that $g$ is flat, 
$g_*L$ is locally free and commutes with restriction to fibers.
We get a meromorphic map
$\phi: X\map \p_S(g_*L)$.  There is thus a normal, bimeromorphic model
$\pi:X'\to X$ such that $\phi\circ \pi:X'\to  \p_S(g_*L)$ is a morphism.

After replacing $X$ by $X'$ and again passing to an open subset of $S$, we may assume that $g$ is flat, 
$g_*L$ is locally free, commutes with restriction to fibers, and
$\phi: X\to \p_S(g_*L)$ is a morphism. 
Let $Y\subset \p_S(g_*L)$ denote its image and
$W\subset X$ the Zariski closed set of points where $\pi:X\to Y$ is not smooth.
Set $Y^\circ:=Y\setminus \phi(W)$ and $X^\circ:=X\setminus \phi^{-1}(\phi(W))$.
The restriction $\phi^\circ: X^\circ\to Y^\circ$ is then smooth and proper.

We assume that $\phi^{-1}(y)$ is a single point for a dense set in $Y$, hence
for a  dense set in $Y^\circ$. Since $\phi^\circ$ is  smooth and proper, it is then an isomorphism.   Thus $\phi$ is bimeromorphic on every irreducible fiber that has a nonempty intersection with $X^\circ$. \qed

\begin{cor} \label{K.vb.Zd.M} Let $g:X\to S$ be a proper morphism of   normal, irreducible  analytic spaces. Then $\operatorname{GT}_S(X)\subset S$ is
\begin{enumerate}
\item either  nowhere dense (in the analytic Zariski topology),
\item or  it contains a dense open subset of $S$, and
$g:X\to S$ is Moishezon.
\end{enumerate}
\end{cor}

Proof. Using resolution, we may assume that $X$ is smooth.
By passing to an open subset of $S$, we may also assume that $S$ and $g$ are smooth.
By \cite{MR2242631}
there is an $m$ (depending only on $\dim X_s$) such that $|mK_{X_s}|$ is very big whenever
$X_s$ is of general type. Thus (\ref{L.vb.Zd.M}) applies to
$L=mK_X$. \qed

\medskip
The following is essentially proved in 
\cite[Thm.1.4]{rao-tsa-1} and \cite[Thm.1.2]{rao-tsa}.

\begin{thm} \label{RT.vb.Zd.M} Let $g:X\to S$ be a smooth, proper morphism of   normal, irreducible  analytic spaces.  
Then $\operatorname{MO}_S(X)\subset S$ is
\begin{enumerate}
\item either  contained in a countable union $\cup_i Z_i$, where
$Z_i\subsetneq S$ are  Zariski closed,
\item or  $\operatorname{MO}_S(X)$  contains a dense, open subset of $S$.
\end{enumerate}
Furthermore, if  $ R^2g_*\o_X$ is torsion free then {\rm (2)} can be replaced by
\begin{enumerate}\setcounter{enumi}{2}
\item $\operatorname{MO}_S(X)=S$ and $g$ is locally Moishezon.
\end{enumerate}
\end{thm}

{\it Remark \ref{RT.vb.Zd.M}.4.}  A positive answer to
Conjecture~\ref{lim.of.M.M.conj} for smooth morphisms  would imply  that
in  fact $\operatorname{MO}_S(X)=S$ always holds  in case  (\ref{RT.vb.Zd.M}.2); see (\ref{RT.vb.Zd.M.cor}). 
\medskip

Proof. Assume first that  $ R^2g_*\o_X$ is torsion free.

As in \cite[3.15]{rao-tsa}, the push-forward of the exponential sequence 
$$
0\to \z_X\to \o_X\stackrel{exp}{\longrightarrow} \o_X^\times\to 1
$$
gives
$$
R^1g_*\o_X^\times\to  R^2g_*\z_X  \stackrel{e_2}{\longrightarrow} R^2g_*\o_X.
$$
We may  pass to the universal cover of $S$ and assume that
$ R^2g_*\z_X$ is a trivial  $H^2(X_s, \z)$-bundle. 

Let $\{\ell_i\}$  be those global sections of $ R^2g_*\z_X$ such that
$e_2(\ell_i)\in H^0(S, R^2g_*\o_X)$ is    identically 0, and
 $\{\ell'_j\}$ the other  global sections.
The $\ell_i$ then lift back to global sections of
$R^1g_*\o_X^\times $, hence to line bundles $L_i$ on $X$.

If there is an $L_i$ such that
$\operatorname{VB}_S(L_i)$ contains a dense open subset of $S$, 
then $X\to S$ is Moishezon by (\ref{L.vb.Zd.M}) and we are done.
Otherwise we claim that 
$$
\operatorname{MO}_S(X)\subset \cup_i \operatorname{VB}_S(L_i) \bigcup 
\cup_j \bigl(e_2(\ell'_j)=0\bigr).
\eqno{(\ref{RT.vb.Zd.M}.5)}
$$
To see this assume that $s\notin \cup_j \bigl(e_2(\ell'_j)=0\bigr)$.
Then every line bundle on $X_s$ is numerically equivalent to some
$L_i|_{X_s}$. Since being big is preserved by numerical equivalence,
we see that $X_s$ has a big line bundle  $\Leftrightarrow$ 
$L_i|_{X_s}$ is big for some $i$  $\Leftrightarrow$ 
$L_i|_{X_s}$ is very big for some $i$. 
This completes the case when $ R^2g_*\o_X$ is torsion free.

In general, the torsion subsheaf of $ R^2g_*\o_X$ is supported on a Zariski closed, proper subset, hence (\ref{RT.vb.Zd.M}.2) gives  that if (\ref{RT.vb.Zd.M}.1) does not hold then $\operatorname{MO}_S(X)$ contains a Zariski dense open subset of $S$. 
\qed

\begin{cor} \label{RT.vb.Zd.M.cor} Let $g:X\to S$ be a smooth, proper morphism of   normal, irreducible  analytic spaces whose fibers are Moishezon.
Then  $g$ is locally Moishezon.
\end{cor}

Proof. If $X_s$ is Moishezon, then  Hodge theory 
(\ref{basic.thms.say}.6) tells us that $H^i(X_s, \c)\to H^i(X_s, \o_{X_s})$ is surjective for every $i$. Thus
 $R^2g_*\o_X$ is locally free by (\ref{dbj.lem.lem}), hence
(\ref{RT.vb.Zd.M}.3) applies. \qed
\medskip

 There are many complex manifolds for which Hodge decomposition holds; these are called  
{\it cohomologically K\"ahler} manifolds or $\partial\bar\partial$-manifolds.
We also get the following variant.

\begin{cor} \label{RT.vb.Zd.M.cor.2} Let $g:X\to S$ be a smooth, proper morphism of   normal, irreducible  analytic spaces. Assume that 
$\operatorname{MO}_S(X)$  contains a dense, open subset of $S$ and
all fibers are cohomologically K\"ahler.
Then  $g$ is locally Moishezon. \qed
\end{cor}

We have used the following result of \cite{dub-jar}; see also
\cite[3.13]{MR946250} and  \cite[2.64]{k-modbook}. 

\begin{thm}\label{dbj.lem.lem}  Let  $g:X\to S$ be a smooth, proper morphism of    analytic spaces. Assume that $H^i(X_s, \c)\to H^i(X_s, \o_{X_s})$ is surjective for every $i$  for some $s\in S$. Then  $ R^ig_*\o_X$ is locally free in a neighborhood of $s$ for every $i$. \qed
\end{thm}

\medskip
(Note that the proof in  \cite{dub-jar} works by descending induction on $i$, so although we are interested in  the $i=2$ case, we need the surjectivity of 
$H^i(X_s, \c)\to H^i(X_s, \o_{X_s})$ for every $i\geq 2$.) 
\medskip

\begin{exmp} \label{bad.fams.exmp}

(\ref{bad.fams.exmp}.1) Let $X\to D^{20}$  be a universal family of K3 surfaces.
A smooth, compact surface is Moishezon iff it is projective.
The projective fibers of  $X\to D^{20}$ correspond to a countable union of hypersurfaces
$H_{2g}\subset D^{20}$.

(\ref{bad.fams.exmp}.2) Let $E\subset \p^2$ be a smooth cubic. Fix $m\geq 10$ and let  $X\to D$ be the universal family of surfaces obtained by blowing up $m$ distinct points $p_i\in E$, and then contracting the birational transform of $E$.  (So $D$ is open in $E^m$.) Such a surface is projective iff there are positive $n_i$ such that
$\sum_i n_i[p_i]\sim  nH$  where $H$ is the line class on $\p^2$ and $n=\frac13 \sum_i n_i$. 

Here $X\to D$ is Moishezon and the projective fibers 
correspond to a countable union of hypersurfaces
$H_i\subset D$.  All fibers have log canonical singularities.

\end{exmp}

\section{1-parameter families}\label{sec.4}

\begin{defn} \label{fib.w.bim.defn} Let $g_i:X^i\to S$  be a  proper morphisms. 
A bimeromorphic map  $\phi:X^1\map X^2$ is 
{\it fiberwise bimeromorphic} if $\phi$ induces a 
bimeromorphic map  $\phi_s: X^1_s\map X^2_s$ for every $s\in S$. 

If $X^1, X^2$ are fiberwise bimeromorphic then 
$X^1_s, X^2_s$ are bimeromorphic to each other for every $s\in S$, but
this is only a sufficient condition in general.

We study whether  a flat, proper, Moishezon  morphism  $g:X\to \dd$
is  fiberwise bimeromorphic to a flat, projective morphism  $g^{\rm p}: X^{\rm p}\to \dd$.
The next examples suggest that the answer is
\begin{itemize}
\item    negative if $g$ is very singular, 
\item  positive if $g$ is mildly singular, and
\item  even if $g$ is smooth, 
$g^{\rm p}$  usually can not be chosen smooth.
\end{itemize}
\end{defn}

\begin{exmp} \label{np.m.exmp.2}
 Let  $g:X\to \dd$ be a smooth, projective morphism.
Assume that $\pic(X)\cong \z$ but $\rank \pic(X_0)\geq 2$.

Let $Z\subset X_0$ be a smooth, ample divisor whose class is not in the image of $\pic(X)\to \pic(X_0)$. Blow up $Z$ to get
$g':X'\to \dd$. Here $X'_0\cong X_0$ has normal bundle $\o_{X_0}(-Z)$, hence it is contractible. We get  a  non-projective, Moishezon morphism  $h:Y\to \dd$.

{\it Conjecture \ref{np.m.exmp.2}.1.}  In most cases, $h:Y\to \dd$ is not fiberwise birational to a  flat,   projective morphism.
\end{exmp}

The next result is a positive answer to the log terminal case of
(\ref{proj.mod.conj.c/lt}), provided
$X_0$ is not uniruled. 
See (\ref{can.main.question.1p.say}) for a discussion of the canonical case.

\begin{thm}\label{proj.mod.thm.lt}
 Let $g:X\to \dd$ be a flat, proper,   Moishezon morphism. Assume that  
\begin{enumerate}
\item  $X_0$ has log terminal  singularities and 
\item $X_0$ is not uniruled.
\end{enumerate}
 Then  $g$ is  
 fiberwise birational to a  flat,  projective morphism
  $g^{\rm p}: X^{\rm p}\to \dd$ (possibly over a  smaller disc) such that
\begin{enumerate}\setcounter{enumi}{2}
\item  $X_0^{\rm p}$ has log terminal  singularities, 
\item $X_s^{\rm p}$ is not uniruled and has terminal singularities  for $s\neq 0$, and
\item $K_{X^{\rm p}}$ is $\q$-Cartier.
\end{enumerate}
\end{thm}

{\it Remark \ref{proj.mod.thm.lt}.6.} Conjecturally we can also achieve that
$K_{X^{\rm p}}$ is  relatively nef. The main obstacle is that (algebraic) minimal models are currently  known to exist only in the general type case.
\medskip

\begin{say}[Proof of (\ref{proj.mod.thm.lt})] \label{proj.mod.thm.lt.pf}  The basic plan is similar to the proof of properness of the KSB moduli space;
see \cite[Sec.5]{ksb} or  \cite[Sec.2.5]{k-modbook}.

 We take a resolution of singularities
$Y\to X$ such that $Y\to\dd$ is projective, and then
take a relative minimal model of   $Y\to\dd$. We hope that it gives what we want.
There are, however,  several  obstacles.  Next we discuss these, and their solutions, but for all technical details we refer to later sections.

(\ref{proj.mod.thm.lt.pf}.1) 
We need to control the singularities of $X$. 
First (\ref{can.modif.2.cor}) reduces us to the case when
$K_X$ is $\q$-Cartier. We assume this from now on. 
Then (\ref{inv.adj.k.thm}) implies that the pair
$(X, X_0)$ is plt.

(\ref{proj.mod.thm.lt.pf}.2) After a base change   $z\mapsto z^r$ we get
$g^r:X^r\to \dd$. For suitable $r$,  there is 
 a semi-stable, projective resolution  
$h:Y\to \dd$; we may also choose it to be equivariant for the
action of the cyclic group $G\cong \z_r$. All subsequent steps will be $G$-equivariant.  We denote by $X_0^Y$ the birational transform of $X_0$ and by $E_i$ the other irreducible components of $Y_0$.

(\ref{proj.mod.thm.lt.pf}.3) We claim that $Y_s$ is not uniruled for $s\neq 0$. Indeed, for smooth families being uniruled is a deformation invariant property, and by Matsusaka's theorem \cite[IV.1.7]{rc-book}, we would get that $X^Y_0$ is uniruled. Thus $K_{Y_s}$ is pseudo-effective by \cite{d-etal}.

(\ref{proj.mod.thm.lt.pf}.4) The required relative  minimal model theorem is known only when  the general fibers are of general type. To achieve this,
let $H$ be an  ample, $G$-equivariant divisor such that
$Y_0+H$ is snc. For $\epsilon>0$ we get  a pair 
$(Y, \epsilon H)$ whose  general fibers $(Y_s, \epsilon H_s)$ are of   log general type since $K_{Y_s}$ is pseudo-effective. For such algebraic families,  relative minimal models are known to exist \cite{bchm}. We also know that $(Y, Y_0+\epsilon H)$ is dlt for $0<\epsilon\ll 1$.

(\ref{proj.mod.thm.lt.pf}.5) Although our family is not algebraic,
\cite{k-nx} treats the relative MMP 
for semi-stable, projective morphisms to a disc.  The precise results are recalled in (\ref{m.c.models.proj.prop}).  Thus we get a relative minimal model
$$
\phi: (Y, \epsilon H)\map (Y^{\rm m}, \epsilon H^{\rm m}),
$$
and $(Y^{\rm m},Y^{\rm m}_0+\epsilon H^{\rm m})$ is dlt. Here $H^{\rm m}$ is $\q$-Cartier for general choice of $\epsilon$  by 
\cite[Lem.1.5.1]{MR3380944}, thus  $(Y^{\rm m},Y^{\rm m}_0)$ is also dlt.
 
{\it Remark.} We have a choice here whether to take the minimal or the canonical model. The minimal model has milder singularities, but it is not unique.
Conjecturally, the canonical model $Y^{\rm c} $ is independent of $0<\epsilon\ll 1$, but this is known only in dimensions $\leq 3$.

(\ref{proj.mod.thm.lt.pf}.6)  We claim that $\phi$ contracts all the $E_i$. 
Since $(X^r, X_0)$ is plt, all the $E_i$ have discrepancy $>-1$. 
Thus the $E_i$ are contained in the restricted, relative  base locus  
of $K_Y+Y_0$ by  (\ref{rest.bs.say}.2).  For $\epsilon$ small enough, the
$E_i$ are also contained in the  restricted, relative  base locus  
of $K_Y+Y_0+\epsilon H$ by  (\ref{rest.bs.say}.1).  Thus any MMP contracts the $E_i$.  On the other hand, $X_0^Y$ can not be contracted, so
$X\map Y^{\rm m}$ is fiberwise birational.

(\ref{proj.mod.thm.lt.pf}.7) Note that $h$ is smooth away from $Y_0$, thus $ (Y_s,\epsilon H_s)$
is terminal for  $s\neq 0$ and  $0\leq \epsilon\ll 1$. Since $H_s$ is ample, we do not contract it, so $(Y^{\rm m}_s, \epsilon H^{\rm m}_s)$  is still terminal. Hence so is   $Y^{\rm m}_s$, giving (4).

(\ref{proj.mod.thm.lt.pf}.8)  As we noted, 
$(Y^{\rm m},Y^{\rm m}_0)$ is dlt, hence plt since
$Y^{\rm m}_0$ is irreducible. Thus $Y^{\rm m}_0$ is log terminal by the easy direction of 
(\ref{inv.adj.k.thm}). 
\qed
\end{say}

\medskip

The following results were also used in the proof of (\ref{proj.mod.thm.lt}).

\begin{prop}[Inversion of adjunction I] \label{inv.adj.k.thm}
 Let $X$ be a normal, complex analytic space,
$X_0\subset X$ a Cartier divisor and $\Delta$ an effective $\r$-divisor
such that $K_X+\Delta$ is $\r$-Cartier. 
Then $(X, X_0+\Delta)$ is plt in a neighborhood of $X_0$ iff
$(X_0, \Delta|_{X_0})$ is klt. 
\end{prop}

Proof. The proof  given in
\cite[Sec.17]{k-etal} or \cite[Sec.5.4]{km-book} applies with minor changes,
using the  complex analytic   vanishing theorems proved in \cite{takegoshi1985} and \cite{MR946250}. \qed

\begin{say}[Divisorial restricted base locus]\label{rest.bs.say} The basic theory is in \cite[1.12--21]{MR2530849} and  an extension to the non-projective case is outlined in \cite[Sec.5]{fkl-vol}.

Let $g:X\to S$ be a proper, Moishezon  morphism, $X$ normal.
The {\it base locus} of a Weil divisor $F$ is 
$$
B(F):=\supp\coker\bigl[g^*g_*\o_X(F)\to \o_X(F)\bigr].
$$
Its divisorial part is denoted by $B^{\rm div}(F)$; we think of it as a Weil divisor.

Let $D$ be an $\r$-divisor on $X$.
Its {\it stable divisorial base locus} is the $\r$-divisor 
$$
{\mathbf B}^{\rm div}(D):=\lim_{m\to\infty}  \tfrac1{m}B^{\rm div}(\rdown{mD}),
$$ and 
 its {\it  restricted  divisorial base locus}
is 
$$
{\mathbf B}^{\rm div}_{-}(D):=\sup_A {\mathbf B}^{\rm div}(D+A),
$$
where $A$ runs through all big $\r$-divisors on $X$ that satisfy
${\mathbf B}^{\rm div}(A)=\emptyset$. This could be an infinite linear combination of prime divisors.

An important observation of \cite{fkl-vol} is that all the projective theorems on the  {\em divisorial} restricted base locus carry over to 
 proper schemes and Moishezon varieties.
We need 2 properties:

(\ref{rest.bs.say}.1) Let $X\to S$ be a proper, Moishezon  morphism,   $D$ an  $\r$-divisor on $X$, and $A$ a big $\r$-divisor on $X$ such that
 ${\mathbf B}^{\rm div}(A)=\emptyset$.
  Then,  
for every prime divisor $F\subset X$,
$$
\coeff_F {\mathbf B}^{\rm div}_{-}(D)=\lim_{\epsilon\to 0} \coeff_F {\mathbf B}^{\rm div}_{-}(D+\epsilon A).
$$

(\ref{rest.bs.say}.2) Let $X_i\to S$ be proper, Moishezon  morphisms,  $h:X_1\to X_2$ a proper, bimeromorhic morphism, $D_2$ a pseudo-effective, $\r$-Cartier divisor on $X_2$, 
and $E$ an effective, $h$-exceptional divisor.  Then
$$
{\mathbf B}^{\rm div}_{-}(E+h^*D_2)\geq E.
$$
\end{say}

\begin{say}[Canonical case of Conjecture~\ref{proj.mod.conj.c/lt}] \label{can.main.question.1p.say}

An  argument similar to (\ref{proj.mod.thm.lt.pf}) should prove  the canonical case, but there are 3 difficulties.

The reduction to the case when $K_X$ is $\q$-Cartier again follows from
(\ref{can.modif.2.cor}). Then we need to show that the pair $(X, X_0)$ is canonical.  This is  proved (though not stated) in 
\cite[5.2]{nak-book}.
This is also a special case of  the general inversion of adjunction; a quite roundabout proof for Moishezon morphisms is given in (\ref{inv.adj.thm}).

 In (\ref{proj.mod.thm.lt.pf}) next  we run the MMP for     $K_Y+\epsilon H$,
which is the same as MMP for     $K_Y+Y_0+\epsilon H$  since $Y_0$ is numerically relatively trivial.

In the canonical case we would need to run the MMP for
$K_Y+X^Y_0+\eta \sum E_i+\epsilon H$, where we choose
$\epsilon, \eta$ small, positive. The arguments of \cite{k-nx} do not cover this case, but I expect that a method similar to \cite{k-nx} would prove this.

If $X_0$ is of general type, then the canonical model of 
$(Y, X^Y_0)$ gives what we want. 

In general, arguing as in (\ref{m.c.models.proj.prop})  we should get a minimal model   $g^{\rm m}:(Y^{\rm m},Y^{\rm m}_0+ \epsilon H^{\rm m})\to \dd$. Here $\eta \sum E_i^{\rm m} $ is omitted since the 
$E_i$ get contracted. General theory tells us that 
 $$
\discrep(Y^{\rm m}_0)\geq \discrep(Y^{\rm m}_0,\diff_{Y^{\rm m}_0}\epsilon H^{\rm m})\geq -\epsilon.
$$
We can choose $\epsilon$ arbitrarily small, but $Y^{\rm m}_0 $ may depend on
$\epsilon$, so we can not just take a limit as $\epsilon\to 0$. This is a problem
that appears even if we start with a projective,  algebraic family.

At this point we could appeal to one of the ACC conjectures 
(\ref{gap.conj.say}) which says that, for $\epsilon$ small enough,  we must have
$\discrep(Y^{\rm m}_0)\geq 0$. That is, $Y^{\rm m}_0$ is canonical. 

The necessary result is known in dimensions $\leq 3$, but it  is likely to be quite difficult in general.  So an alternate approach to our situation would be better.
\end{say}

\begin{say}[A gap conjecture]\label{gap.conj.say} 
The following is a special case of
\cite[Conj.4.2]{MR1420223}

\begin{enumerate}
\item For every $n\geq 1$ there is an $\epsilon(n)>0$ such that
if $X$ is an $n$-dimensional variety and
 $\discrep X>-\epsilon(n)$, then in fact 
 $\discrep X\geq 0$ (that is, $X$ has canonical singularities). 
\end{enumerate}

In dimension 2 this can be read off from the classifiation of log terminal singularities (these have $\discrep X>-1$). We get the optimal value
$\epsilon(2)=\frac1{3}$ and equality holds for $\c^2/\frac1{3}(1,1)$. 

The 3-dimensional case is much harder; see
\cite{jia-gap}. The optimal value  is $\epsilon(3)=\frac1{13}$ and
the extremal case is the cyclic
quotient singularity  $\c^3/\frac1{13}(3, 4, 5)$. 

 Special cases (in all dimensions) are proved in \cite{nakamura2016}.
\end{say}

\begin{rem}\label{def.i.pg.say} 
The deformation invariance of plurigenera for smooth, proper morphisms
with  Moishezon fibers is proved in \cite[Thm.1.2]{rao-tsa}. 

The canonical case of Conjecture~\ref{proj.mod.conj.c/lt}
would show that the projective case implies the  Moishezon case. 
However,  the hard part in \cite{rao-tsa} is to show that $g$ is Moishezon,
so using Conjecture~\ref{proj.mod.conj.c/lt} would only yield a longer proof. 
\end{rem}

\section{Approximating Moishezon  morphisms}\label{sec.5}

We discuss 2 ways of approximating a projective (resp.\ Moishezon) morphism
$g:X\to \dd$ by morphisms between projective (resp.\ Moishezon) varieties.
This allows us to prove some results  for Moishezon morphisms
$g:X\to \dd$.

\begin{say}[Algebraic approximation of projective morphisms]\label{a.approx.m.1p.lem}
Let $g:Y\to \dd$ be a projective morphism with relatively ample line bundle $L$.  For later purposes we also  specify a finite set of
relative Cartier divisors $E^i\subset Y$.

Then  $(Y_0, L_0:=L|_{Y_0},  E^i_0:=E^i|_{Y_0})$ is a projective, polarized scheme marked with effective Cartier divisors.  (For now $Y_0$  can be even nonreduced.)

$(Y_0, L_0,  E^i_0)$  has a universal deformation space
$G_S: ({\mathbf Y}_S, {\mathbf L}_S,  {\mathbf E}^i_S)   \to S$,
where $S, {\mathbf Y}_S$ are quasi-projective schemes, $G_S$ is flat and projective, ${\mathbf L}_S$ is $G_S$-ample  and the ${\mathbf E}^i_S$ are relative Cartier divisors.

The original family gives a holomorphic   $\phi_S:\dd\to S$. 
Next we replace $S$ first by the Zariski closure of $\phi_S(\dd)$
and then by its resolution. We obtain the following data.
\begin{enumerate}
 \item A smooth $\c$-variety $B$,
\item a flat, projective morphism $G:({\mathbf Y}, {\mathbf L},  {\mathbf E}^i)   \to B$, where ${\mathbf L}$ is $G$-ample,  the ${\mathbf E}^i$ are relative Cartier divisors, and
\item a holomorphic  map  $\phi:\dd\to B(\c)$,
\end{enumerate}
such that,
\begin{enumerate}\setcounter{enumi}{3}
\item   $\bigl((Y,L,  E^i)\to \dd\bigr)\cong \bigl(\phi^*({\mathbf Y}, {\mathbf L},  {\mathbf E}^i)   \to \dd\bigr)$ and
\item  $\phi(\dd)$ is smooth and Zariski dense in $B$.
\end{enumerate}
We call  $G:({\mathbf Y}, {\mathbf L}, {\mathbf E}^i)   \to B$ an {\it algebraic envelope} of  $g:(Y,L,  E^i)\to \dd$.

Note that  we have no control over the dimension of $B$.
However, if $Y$ is smooth, then so is ${\mathbf Y} $.

Since $B$ is smooth, the holomorphic curve $\phi(\dd)\subset B$ can be approximated by algebraic curves to any order. Thus, for any fixed $m>0$  we get 
\begin{enumerate}\setcounter{enumi}{5}
 \item a smooth, pointed algebraic curve  $(c, C)$,
\item a flat, projective morphism $(Y_C, L_C,  E^i_C)   \to C$, and
\item an isomorphism 
  $(Y,L,  E^i)_m\cong (Y_C, L_C, E^i_C)_m$,
\end{enumerate}
where the subscript $m$ denotes the $m$th order infinitesimal neighborhood of
the central fibers.

If $m=0$ then we only get that the central fibers  $Y_0$ and $(Y_C)_0$ are isomorphic. The case $m=1$ carries much more information: the smoothness of the total space along the central fiber and the normal bundles of the irreducible components of the central fiber are also preserved.
\end{say}

As in \cite{k-nx}, algebraic  envelopes can be used to
show that MMP works for projective morphism over Riemann surfaces.

\begin{prop} \label{m.c.models.proj.prop}
Let $g:(Y, \Delta)\to \dd$ be projective, locally stable with irreducible, normal general fibers.
Assume that $K_Y+\Delta$ is big on almost every fiber. Then 
\begin{enumerate}
\item there is a relative minimal model $g^{\rm m}:(Y^{\rm m}, \Delta^{\rm m})\to \dd$, and
\item the relative canonical model $g^{\rm c}:(Y^{\rm c}, \Delta^{\rm c})\to \dd$ exists.
\end{enumerate}
Assume in addition that $(Y, Y_0+\Delta)$ is dlt and there is an irreducible divisor $Y^*_0\subset Y_0$ 
such that $Y_0\setminus Y_0^*$ is contained in the stable, relative  base locus of   $K_Y+\Delta$. Then
\begin{enumerate}\setcounter{enumi}{2}
\item $Y_0^*\map Y^{\rm m}_0\to Y^{\rm c}_0 $ are birational, 
\item $(Y^{\rm m}, Y^{\rm m}_0+\Delta^{\rm m})$ and $(Y^{\rm c}, Y^{\rm c}_0+\Delta^{\rm c})$ are plt, and
\item $(Y^{\rm m}_0,\diff_{Y^{\rm m}_0}\Delta^{\rm m})$ and $(Y^{\rm c}_0,\diff_{Y^{\rm c}_0}\Delta^{\rm c})$ are klt.
\item If $\Delta$ is $\r$-Cartier then $(Y^{\rm m}, Y^{\rm m}_0)$ is plt and
$Y^{\rm m}_0$ is log terminal.
\item If the coefficients in $\Delta$ are sufficiently general,   then $(Y^{\rm c}, Y^{\rm c}_0)$ is plt and
$Y^{\rm c}_0$ is log terminal.
\end{enumerate}
\end{prop}

Proof. Claims (1--2) are basically proved in  \cite{k-nx}.
Unfortunately, the main result \cite[Thm.2]{k-nx} is formulated to apply to the Calabi-Yau case. However, \cite[Props.8--14]{k-nx} contain a complete proof, though not a clear statement.

Any MMP $Y\map Y^{\rm m}$ contracts the stable base locus of   $K_Y+\Delta$, thus
$Y^{\rm m}_0 $ is irreducible and $Y_0^*\map Y^{\rm m}_0$ is thus bimeromorphic. 
Also $(Y^{\rm m}, Y^{\rm m}_0+\Delta^{\rm c})$ is dlt, hence plt  since 
$Y^{\rm m}_0 $ is irreducible.  Since $Y^{\rm m}\to Y^{\rm c}$ does not contract
$Y^{\rm m}_0 $, we see that $(Y^{\rm c}, Y^{\rm c}_0+\Delta^{\rm c})$ is also plt.
This is (4), and (5) follows by the easy direction of   adjunction
\cite[4.8]{kk-singbook}. 

If $\Delta$ is $\q$-Cartier then so is $\Delta^{\rm m}$, hence (6) follows from \cite[2.27]{km-book}.  A similar argument works for (7) using \cite[Lem.1.5.1]{MR3380944}. 
\qed

\begin{say}[Algebraic approximation of Moishezon morphisms]\label{m.approx.m.1p.lem}
Let $f:X\to \dd$ be a proper, Moishezon morphism and
$h:Y\to X$ a proper morphism such that $g:=f\circ h:Y\to \dd$ is projective with relatively ample line bundle $L$.   Also choose  relative Cartier divisors $E^i$
on $Y$. Assume also that $X_0$ is seminormal (though this is probably ultimately not necessary). 

We can apply (\ref{a.approx.m.1p.lem}.1--5) to get an algebraic envelope 
$G:({\mathbf Y}, {\mathbf L},   {\mathbf E}^i)   \to B$. By \cite{artin}, after an \'etale base change, we may assume that   $h_0:Y_0\to X_0$ extends to 
 $H: {\mathbf Y}\to {\mathbf X}$ where $F:{\mathbf X}\to B$ is an algebraic space.

{\it Comment \ref{m.approx.m.1p.lem}.1.} General extension theory, as in \cite{artin, MR0304703},  tells us only that we have  
$$
 H_0: {\mathbf Y}_0=Y_0\stackrel{h_0}{\longrightarrow}  X_0
\stackrel{\tau}{\longrightarrow}{\mathbf X}_0,
$$
where $\tau$ is a finite homeomorphism.
Then we use  that the functor of simultaneous seminormalizations is formally representable. The projective case is discussed in \cite{k-hh-2}; see
\cite[9.61]{k-modbook} for algebraic spaces. 
We assumed that we have a
simultaneous seminormalization over the completion of $\phi(\dd)$,
which is Zariski dense. Thus ${\mathbf X}\to B$ has seminormal fibers, hence
$X_0\cong {\mathbf X}_0$, as claimed.

Assume next that  fibers of $f$ over $\dd^\circ$ satisfy a property ${\mathcal P}$ that is Zariski open in families (for example smooth, normal or reduced). Then general fibers of $F$ also satisfy  ${\mathcal P}$.
 As before,  $\phi(\dd)\subset B$ can be approximated by algebraic curves to any order. Thus, for any fixed $m>0$  we get
\begin{enumerate}\setcounter{enumi}{1}
 \item a smooth, pointed algebraic curve  $(c, C)$,
\item  morphisms $h_c:(Y_C, L_C,  E^i_C) \to X_C   \to C$, where
\begin{enumerate} 
 \item $g_C:(Y_C, L_C)    \to C$ is flat, projective,
\item $f_C: X_C  \to C$ is a flat algebraic space, 
\item  general fibers of $f_C$ satisfy  ${\mathcal P}$, and
\end{enumerate}
\item an isomorphism 
  $\bigl((Y,L,  E^i)\to X\bigr)_m\cong \bigl((Y_C, L_C,  E^i_C)\to X_C\bigr)_m$,
\end{enumerate}
where the subscript $m$ denotes the $m$th order infinitesimal neighborhood of
the central fibers.

\end{say}

\begin{cor} \label{can.modif.1.cor} 
Let $f:X\to \dd$ be a flat, proper, Moishezon morphism, $X$ normal.
Assume that it has a resolution 
$h:Y\to X$ where  $g:=f\circ h:Y\to \dd$ is projective and $Y_0$ a reduced, snc divisor. Then $X$ has a canonical modification $\pi: X^{\rm c}\to X$.
(That is, $X^{\rm c}$ has canonical singularities and $K_{ X^{\rm c}}$ is $\pi$-ample.)
\end{cor}

Proof. Let   $H: {\mathbf Y}\to {\mathbf X}$ and  $F:{\mathbf X}\to B$ be an algebraic envelope as in  (\ref{m.approx.m.1p.lem}). 

Note that canonical modifications are unique and commute with \'etale morhisms. They exist for quasi-projective varieties over $\c$ by \cite{bchm}, hence every algebraic space of finite type over $\c$ has a canonical modification.

Let
$\Pi:{\mathbf X}^{\rm c}\to {\mathbf X}$ denote  the canonical modification of ${\mathbf X}$.
Since  ${\mathbf Y}\to B$ is locally stable, so is 
${\mathbf X}^{\rm c}\to B$; cf.\ \cite[Sec.4.8]{k-modbook}. 

By pull-back we get  a locally stable morphism $\pi: X^{\rm c}\to X\to \dd$
whose general fibers are canonical. Since $(X^{\rm c}, X^{\rm c}_0)$ is lc and
$X^{\rm c}_0$ is a  Cartier divisor, we see that $X^{\rm c}$ has canonical singularities.\qed

\medskip

The following extends \cite[Sec.3]{ksb} to Moishezon morphisms, see also \cite[Sec.5.5]{k-modbook}.

\begin{cor} \label{can.modif.2.cor} 
Let $f:X\to \dd$ be a flat, proper, Moishezon morphism.
Assume that  $X_0$ is log terminal. 
Then $X$ has a canonical modification $\pi: X^{\rm c}\to X$,
$X^{\rm c}_0$ is log terminal and $\pi$ is fiberwise birational.
\end{cor}

Proof.  After a ramified base change $\tilde \dd\to \dd$ with group $G:=\z/r$, we can apply (\ref{can.modif.1.cor}) to  $\tilde X\to \tilde \dd$ to get
$\tilde\pi: \tilde X^{\rm c}\to \tilde  X$.  

As in \cite[5.32]{k-modbook} we get that $\tilde X^{\rm c}_0$ is log terminal
and  $\tilde X^{\rm c}_0\to \tilde X_0$ is birational. 
Set $X^{\rm c}:=\tilde X^{\rm c}/G$. 

The base change group acts trivially on the central fiber  $\tilde X^{\rm c}_0$, hence  $X^{\rm c}_0\cong \tilde X^{\rm c}_0$ is also log terminal. 
Finally  the pair $\bigl(\tilde X^{\rm c}, \tilde X^{\rm c}_0\bigr)$ is log canonical, hence so is $\bigl(X^{\rm c}, X^{\rm c}_0\bigr)$ by \cite[5.20]{km-book}.
Thus $X^{\rm c}$ is  canonical. 
\qed

\section{Inversion of adjunction}\label{sec.7}

The proof of the general
inversion of adjunction theorem  given in \cite[4.9]{kk-singbook} 
relies on MMP, which is not known for projective morphisms over an analytic base. (See \cite{MR946250} for the first steps and \cite{k-defgt} for some special cases.)

We go around this for Moishezon morphisms using approximations.

\begin{thm} \label{inv.adj.thm}
 Let $g:X\to \dd$ be a flat, proper, Moishezon morphism and $\Delta$ an effective $\q$-divisor on $X$. Assume that $K_X+\Delta$ is $\q$-Cartier.
Then 
$$
\discrep(X, X_0+\Delta_0)= \totaldiscrep(X_0, \Delta_0),
$$
where on the left   we use only those exceptional divisors $E$ over $X$ whose centers on $X$ have nonempty intersection with  $X_0$. 
\end{thm}

Note that here the $\leq $ part is easy \cite[4.8]{kk-singbook}. The known proofs of the $\geq $ part  use MMP, and the cases settled in \cite{k-nx} do not seem  enough.

We start the proof of (\ref{inv.adj.thm}) with a discussion on snc divisors and then 
with a general result which says that discrepancies can be computed from the 1st order neighborhood of the exceptional set.

\begin{say}[Simple normal crossing divisors]\label{simp.norm.say}
It would be convenient  to recognize simple normal crossing divisors
(abbreviated as {\it snc})
from an infinitesimal neighborhood of the special fiber.  At first sight, this seems impossible. Consider for example the family
$$
g:\bigl(\c^3, D:=(xy=z^{m+1})\bigr)\to \c_z,
$$
where $D$ is not an snc divisor.
The $m$th order infinitesimal neighborhood of the special fiber is defined by
$z^{m+1}=0$, hence isomorphic to the $m$th order  neighborhood
of the snc family
$$
g:\bigl(\c^3, B:=(xy=0)\bigr)\to \c_z.
$$
There is also the added problem that snc in the Euclidean topology is not the same as snc in the Zariski topology. (For example,  $(y^2=x^2+x^3$ is snc in the Euclidean topology but not in the Zariski topology.)

We can, however, solve both problems by a simple bookkeeping convention.

Let  $M$ be a complex manifold and $\{E_i:i\in I\}$ (reduced) divisors on $M$.
We say that  $(M,  E_i: i\in I)$ is a {\it marked snc pair} if for every $p\in M$ there are 
\begin{enumerate}
\item local analytic coordinates $z_1, \dots, z_n$,  and
\item  an injection $\sigma:\{1,\dots, r\}\into I$  for some $0\leq r\leq m$,
\end{enumerate}
such that
\begin{enumerate}\setcounter{enumi}{2}
\item  $E_{\sigma(i)}=(z_i=0)$ near $p$, and
\item the other $E_j$  do not contain $p$.
\end{enumerate}
With this definition we have the following.

\medskip
{\it Claim \ref{simp.norm.say}.5.} Let 
$E_1,\dots, E_r, E_{r+1}, \dots, E_m$ and $E'_{r+1}, \dots, E'_m$ be reduced divisors on a   complex manifold $M$.  Assume that
\begin{enumerate}
\item[(a)] $(M, E_1+\cdots +E_m)$ is a marked snc pair, and
\item[(b)]  $E_j$ and $E'_j$ have the same restriction on 
$E_1\cup\dots\cup E_r $ for all $j>r$.
\end{enumerate}
Then  $(M, E_1+\cdots +E_r+E'_{r+1}+\cdots +E'_m)$ is also a marked snc pair in a neighborhood of $E_1\cup\dots\cup E_r $. \qed
\end{say}

\begin{lem} \label{discrep.2nd.oredr.lem}
 Let $(X, \Delta=\sum d_jD_j)$ be a normal  analytic pair such that
$K_X+\Delta$ is $\q$-Cartier. Let $B_X\subset X$ be a Cartier divisor.
Let 
$$
p:(Y, \tsum_i B_i+\tsum_j \bar D_j+\tsum_\ell E_\ell)\to (X, B_X+\supp\Delta)
$$ be a log resolution,
where $B:=\sum_i B_i=\red p^{-1}(B_X)$,  $\bar D_j$ is the birational transform of $D_j$ and $E_\ell$ are the  other $p$-exceptional divisors. 

Then the discrepancies $a(*, X, \Delta)$ of all $p$-exceptional divisors 
(whose centers have nonempty intersection with $B_X$) 
can be computed from
\begin{enumerate}
\item    $B_{(2)}:=\spec \o_Y/\o_Y(-2B)$, 
\item   $\bar D_j|_B$ and $E_\ell|_B$. 
\end{enumerate}
\end{lem}

Proof. After  replacing $X$ by a smaller neighbood of $B_X$, we may assume that
$B$ is a deformation retract of $Y$.
In particular, the  centers of all $p$-exceptional divisors have nonempty intersection with $B_X$, and numerical equivalence of divisors is determined by their restriction to $B$. 

The discrepancies  $b_i$ and $e_\ell$ are uniquely determined by the conditions
$$
K_Y+\tsum'_i b_iB_i+\tsum_j d_j\bar D_j+\tsum_\ell e_\ell E_\ell\equiv_p p^*(K_X+\Delta),
\eqno{(\ref{discrep.2nd.oredr.lem}.3)}
$$
where in $\tsum'_i $ we sum over the $p$-exceptional divisor in $B$.
Restricting to $B$ and using adjunction we get
$$
\tsum'_i b_i(B_i|_B)+\tsum_\ell e_\ell (E_\ell|_B)\equiv_p -K_B+(B|_B)+(p|_B)^*(K_X+\Delta).
\eqno{(\ref{discrep.2nd.oredr.lem}.4)}
$$
Note that $B_{(2)}$ determines the $B_i|_B$ and hence $B|_B$. 
Thus the right hand side is known and the $b_i, e_\ell$ are the unique solution to (\ref{discrep.2nd.oredr.lem}.4). \qed

\begin{cor} \label{discrep.2nd.oredr.cor}
Using the notation of (\ref{discrep.2nd.oredr.lem}), assume that there is another pair with a log resolution
$$
p':(Y', \tsum_i B'_i+\tsum_j \bar D'_j+\tsum_\ell E'_\ell+ F')\to (X', B'_{X'}+\supp\Delta')
$$
such that there is an isomorphism
$$
\phi:
\bigl( B'_{(2)}\hookleftarrow B'\stackrel{p'}{\to} B'_{X'}\bigr)
\cong 
\bigl( B_{(2)}\hookleftarrow B\stackrel{p}{\to} B_X\bigr),
$$ 
that sends $  \bar D'_j|_{B'}$ to  $  \bar D_j|_{B}$ and
$E'_\ell|_{B'}$ to   $E_\ell|_{B}$ for every $j, \ell$. Then
\begin{enumerate}
\item corresponding divisors have the same discrepancies, and
\item divisors in $F'$ have discrepancy 0.
\end{enumerate}
\end{cor}

Proof. Note that (\ref{discrep.2nd.oredr.lem}.4) gives us that
$$
\tsum'_i b_i(B'_i|_{B'})+\tsum_\ell e_\ell (E'_\ell|_{B'})+ 0\cdot F'\equiv_{p'} -K_{B'}+(B'|_{B'})+(p'|_{B'})^*(K_{X'}+\Delta').
$$
Since this equation has a unique solution,  $b_i, e_\ell$ give the 
 discrepancies over $X'$. \qed

\medskip

The following example illustrates the role of the divisor $F'$ in 
(\ref{discrep.2nd.oredr.cor}). 

\begin{exmp} \label{new.exc.divs.exmp}
Let $X=(x^2-y^2+z^2=t^4)\subset \c^4$, $B=(t=0)$ and $Y$ the small resolution obtained by blowing up $(x-y=z-t^2=0)$.  (Here $\Delta=0$ and $E$ s empty.) Next set  $X'=(x^2-y^2+z^2=0)\subset \c^4$, $B'=(t=0)$ and $Y'$ the  resolution obtained by blowing up $(x=y=z=0)$. 

The 1st order neigborhoods are isomorphic, but on $Y'$ 
 we have an exceptional divisor $F'$.  
Note that if we replace $t^4$ by $t^{2m+2}$, we have isomorphisms of $m$th order 
infinitesimal neighborhoods as well.

Thus we can not tell whether a singularity is terminal or canonical by looking at $m$th order infinitesimal neighborhoods for some fixed $m$.  
\end{exmp}

\begin{say}[Proof of (\ref{inv.adj.thm})]\label{inv.adj.thm.pf}
Write $\Delta=\sum d_j D_j$. 

Let $h:(Y, \sum \bar D_j)\to (X, \sum D_j)$ be a log resolution of singularities such that $Y\to \dd$ is projective. Let  $\bar D_j$ denote the birational transform of $D_j$, and let  $E_i\subset Y$ be the exceptional divisors that dominate $\dd$.  

 By (\ref{m.approx.m.1p.lem}.2--4) there are a smooth, pointed algebraic curve $(c, C)$, 
a flat, proper morphism of algebraic spaces
$X^{\rm a}\to C$ and a projective resolution  $h^{\rm a}:Y^{\rm a}\to X^{\rm a}$ such that
$$
\bigl( h^{\rm a}:(Y^{\rm a}, \tsum \bar D_j^{\rm a}+\tsum E^{\rm a}_j)\to X^{\rm a}\bigr)_1\cong  \bigl( h:(Y, \tsum D_j+\tsum E_i)\to X\bigr)_1.
\eqno{(\ref{inv.adj.thm.pf}.1)}
$$
Note that $h^{\rm a}(E^{\rm a}_j)\cap X_0^{\rm a}=h(E_j)\cap X_0$, thus
the $E^{\rm a}_j$ are $h^{\rm a}$-exceptional. (As in (\ref{new.exc.divs.exmp}), there may be other 
$h^{\rm a}$-exceptional divisors.)

By (\ref{simp.norm.say}), $(Y^{\rm a}, \tsum \bar D_j^{\rm a}+\tsum E^{\rm a}_j )$ is also an snc pair,
we are thus in the situation of (\ref{discrep.2nd.oredr.cor}). 
Since inversion of adjunction holds for the algebraic pair
$(X^{\rm a}, X_0^{\rm a}+\Delta^{\rm a})$, it also holds for
$(X, X_0+\Delta)$. \qed
\end{say}


\begin{thebibliography}{BCHM10}

\bibitem[AK00]{Abramovich-Karu00}
D.~Abramovich and K.~Karu, \emph{Weak semistable reduction in characteristic
  0}, Invent. Math. \textbf{139} (2000), no.~2, 241--273. \MR{MR1738451
  (2001f:14021)}

\bibitem[Ale15]{MR3380944}
Valery Alexeev, \emph{Moduli of weighted hyperplane arrangements}, Advanced
  Courses in Mathematics. CRM Barcelona, Birkh\"{a}user/Springer, Basel, 2015,
  Edited by Gilberto Bini, Mart\'{\i} Lahoz, Emanuele Macr{\`\i} and Paolo
  Stellari. \MR{3380944}

\bibitem[Art69]{Artin69b}
Michael Artin, \emph{Algebraization of formal moduli. {I}}, Global Analysis
  (Papers in Honor of K. Kodaira), Univ. Tokyo Press, Tokyo, 1969, pp.~21--71.
  \MR{MR0260746 (41 \#5369)}

\bibitem[Art70]{artin}
\bysame, \emph{Algebraization of formal moduli. {II}. {E}xistence of
  modifications}, Ann. of Math. (2) \textbf{91} (1970), 88--135. \MR{0260747
  (41 \#5370)}

\bibitem[Bar75]{MR0399503}
Daniel Barlet, \emph{Espace analytique r\'eduit des cycles analytiques
  complexes compacts d'un espace analytique complexe de dimension finie},
  Fonctions de plusieurs variables complexes, {II} ({S}\'em.
  {F}ran{\cedilla{c}}ois {N}orguet, 1974--1975), Springer, Berlin, 1975,
  pp.~1--158. Lecture Notes in Math., Vol. 482. \MR{0399503}

\bibitem[Bar17]{barlet2017gauduchons}
\bysame, \emph{Gauduchon's form and compactness of the space of divisors},
  2017.

\bibitem[BCHM10]{bchm}
Caucher Birkar, Paolo Cascini, Christopher~D. Hacon, and James
  M\textsuperscript{c}Kernan, \emph{Existence of minimal models for varieties
  of log general type}, J. Amer. Math. Soc. \textbf{23} (2010), no.~2,
  405--468.

\bibitem[BDPP13]{d-etal}
S\'ebastien Boucksom, Jean-Pierre Demailly, Mihai P{\u{a}}un, and Thomas
  Peternell, \emph{The pseudo-effective cone of a compact {K}\"ahler manifold
  and varieties of negative {K}odaira dimension}, J. Algebraic Geom.
  \textbf{22} (2013), no.~2, 201--248.

\bibitem[Cam80]{MR583821}
F.~Campana, \emph{Alg\'ebricit\'e et compacit\'e dans l'espace des cycles d'un
  espace analytique complexe}, Math. Ann. \textbf{251} (1980), no.~1, 7--18.
  \MR{MR583821 (82a:32026)}

\bibitem[Cam81]{MR620706}
\bysame, \emph{R\'eduction alg\'ebrique d'un morphisme faiblement
  {K}\"ahl\'erien propre et applications}, Math. Ann. \textbf{256} (1981),
  no.~2, 157--189. \MR{MR620706 (84i:32014)}

\bibitem[Cam91]{MR1107661}
\bysame, \emph{The class {${ C}$} is not stable by small deformations}, Math.
  Ann. \textbf{290} (1991), no.~1, 19--30. \MR{MR1107661 (93d:32047)}

\bibitem[DJ74]{dub-jar}
Philippe Dubois and Pierre Jarraud, \emph{Une propri\'et\'e de commutation au
  changement de base des images directes sup\'erieures du faisceau structural},
  C. R. Acad. Sci. Paris S\'er. A \textbf{279} (1974), 745--747. \MR{0376678
  (51 \#12853)}

\bibitem[ELM{\etalchar{+}}09]{MR2530849}
Lawrence Ein, Robert Lazarsfeld, Mircea Musta{\cedilla{t}}{\u{a}}, Michael
  Nakamaye, and Mihnea Popa, \emph{Restricted volumes and base loci of linear
  series}, Amer. J. Math. \textbf{131} (2009), no.~3, 607--651. \MR{2530849
  (2010g:14005)}

\bibitem[FKL16]{fkl-vol}
Mihai Fulger, J{\'a}nos Koll{\'a}r, and Brian Lehmann, \emph{Volume and
  {H}ilbert function of {$\Bbb R$}-divisors}, Michigan Math. J. \textbf{65}
  (2016), no.~2, 371--387. \MR{3510912}

\bibitem[Fuj82]{fuj-c}
A.~Fujiki, \emph{On the {D}ouady space of a compact complex space in the
  category {$\mathcal C$}}, Nagoya Math. J. \textbf{85} (1982), 189--211.

\bibitem[Hir75]{MR0393556}
Heisuke Hironaka, \emph{Flattening theorem in complex-analytic geometry}, Amer.
  J. Math. \textbf{97} (1975), 503--547. \MR{MR0393556 (52 \#14365)}

\bibitem[HM06]{MR2242631}
Christopher~D. Hacon and James McKernan, \emph{Boundedness of pluricanonical
  maps of varieties of general type}, Invent. Math. \textbf{166} (2006), no.~1,
  1--25. \MR{MR2242631 (2007e:14022)}

\bibitem[Iit71]{iitaka}
Shigeru Iitaka, \emph{On {$D$}-dimensions of algebraic varieties}, J. Math.
  Soc. Japan \textbf{23} (1971), 356--373. \MR{0285531 (44 \#2749)}

\bibitem[Jia20]{jia-gap}
Chen Jiang, \emph{A gap theorem for minimal log discrepancies of non-canonical
  singularities in dimension three}, 2020.

\bibitem[KM98]{km-book}
J{\'a}nos Koll{\'a}r and Shigefumi Mori, \emph{Birational geometry of algebraic
  varieties}, Cambridge Tracts in Mathematics, vol. 134, Cambridge University
  Press, Cambridge, 1998, With the collaboration of C. H. Clemens and A. Corti,
  Translated from the 1998 Japanese original.

\bibitem[KNX18]{k-nx}
J\'{a}nos Koll\'{a}r, Johannes Nicaise, and Chenyang Xu, \emph{Semi-stable
  extensions over 1-dimensional bases}, Acta Math. Sin. (Engl. Ser.)
  \textbf{34} (2018), no.~1, 103--113. \MR{3735836}

\bibitem[Kol92]{k-etal}
J{\'a}nos Koll{\'a}r (ed.), \emph{Flips and abundance for algebraic
  threefolds}, Soci\'et\'e Math\'ematique de France, 1992, Papers from the
  Second Summer Seminar on Algebraic Geometry, University of Utah, Salt Lake
  City, Utah, August 1991, Ast\'erisque No. 211 (1992).

\bibitem[Kol96]{rc-book}
\bysame, \emph{Rational curves on algebraic varieties}, Ergebnisse der
  Mathematik und ihrer Grenzgebiete. 3. Folge., vol.~32, Springer-Verlag,
  Berlin, 1996.

\bibitem[Kol11]{k-hh-2}
\bysame, \emph{Simultaneous normalization and algebra husks}, Asian J. Math.
  \textbf{15} (2011), no.~3, 437--449. \MR{2838215}

\bibitem[Kol13]{kk-singbook}
\bysame, \emph{Singularities of the minimal model program}, Cambridge Tracts in
  Mathematics, vol. 200, Cambridge University Press, Cambridge, 2013, With the
  collaboration of S{\'a}ndor Kov{\'a}cs.

\bibitem[Kol20]{k-modbook}
\bysame, \emph{Moduli of varieties of general type}, (book in preparation,
  \url{https://web.math.princeton.edu/~kollar/book/modbook20170720.pdf}), 2020.

\bibitem[Kol21]{k-defgt}
\bysame, \emph{Deformations of varieties of general type}, 2021.

\bibitem[KSB88]{ksb}
J{\'a}nos Koll{\'a}r and N.~I. Shepherd-Barron, \emph{Threefolds and
  deformations of surface singularities}, Invent. Math. \textbf{91} (1988),
  no.~2, 299--338.

\bibitem[Loo81]{MR632841}
Eduard Looijenga, \emph{Rational surfaces with an anticanonical cycle}, Ann. of
  Math. (2) \textbf{114} (1981), no.~2, 267--322. \MR{632841 (83j:14030)}

\bibitem[LP92]{leb-poo}
Claude Lebrun and Yat-Sun Poon, \emph{Twistors, {K}\"ahler manifolds, and
  bimeromorphic geometry. {II}}, Journal of the American Mathematical Society
  \textbf{5} (1992), no.~2, 317--325.

\bibitem[Moi66]{Moi-66}
Boris Moishezon, \emph{On $n$-dimensional compact varieties with $n$
  algebraically independent meromorphic functions, {I}, {II} and {III} (in
  {R}ussian)}, Izv. Akad. Nauk SSSR Ser. Mat. \textbf{30} (1966), 133--174,
  345--386, 621--656.

\bibitem[Moi71]{MR0425189}
\bysame, \emph{Algebraic varieties and compact complex spaces}, Actes du
  Congr\`es International des Math\'ematiciens (Nice, 1970), Tome 2,
  Gauthier-Villars, Paris, 1971, pp.~643--648. \MR{MR0425189 (54 \#13146)}

\bibitem[Moi74]{MR0369746}
\bysame, \emph{Modifications of complex varieties and the {C}how lemma},
  Classification of algebraic varieties and compact complex manifolds,
  Springer, Berlin, 1974, pp.~133--139. Lecture Notes in Math., Vol. 412.
  \MR{MR0369746 (51 \#5978)}

\bibitem[MR71]{MR0304703}
A.~Markoe and H.~Rossi, \emph{Families of strongly pseudoconvex manifolds},
  Symposium on {S}everal {C}omplex {V}ariables ({P}ark {C}ity, {U}tah, 1970),
  Springer, Berlin, 1971, pp.~182--207. Lecture Notes in Math., Vol. 184.
  \MR{0304703 (46 \#3835)}

\bibitem[Nak75]{nak-par}
Iku Nakamura, \emph{Complex parallelizable manifolds and their small
  deformations}, J. Diff. Geom. \textbf{10} (1975), 85--112.

\bibitem[Nak87]{MR946250}
Noboru Nakayama, \emph{The lower semicontinuity of the plurigenera of complex
  varieties}, Algebraic geometry, {S}endai, 1985, Adv. Stud. Pure Math.,
  vol.~10, North-Holland, Amsterdam, 1987, pp.~551--590. \MR{946250
  (89h:14028)}

\bibitem[Nak04]{nak-book}
\bysame, \emph{Zariski-decomposition and abundance}, MSJ Memoirs, vol.~14,
  Mathematical Society of Japan, Tokyo, 2004.

\bibitem[Nak16]{nakamura2016}
Yusuke Nakamura, \emph{On minimal log discrepancies on varieties with fixed
  {G}orenstein index}, Michigan Math. J. \textbf{65} (2016), no.~1, 165--187.

\bibitem[Nam02]{nam-02}
Yoshinori Namikawa, \emph{Projectivity criterion of {M}oishezon spaces and
  density of projective symplectic varieties}, Internat. J. Math. \textbf{13}
  (2002), 125--135.

\bibitem[Pop09]{pop09}
Dan Popovici, \emph{Limits of projective and $\partial\bar\partial$-manifolds
  under holomorphic deformations}, 2009.

\bibitem[Pop11]{pop11}
\bysame, \emph{Deformation openness and closedness of various classes of
  compact complex manifolds; examples}, 2011.

\bibitem[Pop19]{pop19}
\bysame, \emph{Adiabatic limit and deformations of complex structures}, 2019.

\bibitem[RT20]{rao-tsa}
Sheng Rao and I-Hsun Tsai, \emph{Invariance of plurigenera and {C}how-type
  lemma}, 2020.

\bibitem[RT21]{rao-tsa-1}
\bysame, \emph{Deformation limit and bimeromorphic embedding of {M}oishezon
  manifolds}, Communications in Contemporary Mathematics \textbf{(to appear)}
  (2021).

\bibitem[Sho96]{MR1420223}
V.~V. Shokurov, \emph{{$3$}-fold log models}, J. Math. Sci. \textbf{81} (1996),
  no.~3, 2667--2699, Algebraic geometry, 4. \MR{MR1420223 (97i:14015)}

\bibitem[{Sta}15]{stacks-project}
The {Stacks Project Authors}, \emph{{S}tacks {P}roject},
  http://stacks.math.columbia.edu, 2015.

\bibitem[Sun80]{sundararaman}
D.~Sundararaman, \emph{Moduli, deformations, and classifications of compact
  complex manifolds}, Pitman Pub. London, 1980.

\bibitem[Tak85]{takegoshi1985}
Kensho Takegoshi, \emph{Relative vanishing theorems in analytic spaces}, Duke
  Math. J. \textbf{52} (1985), no.~1, 273--279.

\bibitem[Uen75]{ueno}
Kenji Ueno, \emph{Classification theory of algebraic varieties and compact
  complex spaces}, Lecture Notes in Mathematics, Vol. 439, Springer-Verlag,
  Berlin, 1975, Notes written in collaboration with P. Cherenack.

\bibitem[Uen83]{ueno-83}
\bysame, \emph{Introduction to the theory of compact complex spaces in the
  class {$\mathcal C$}}, Algebraic Varieties and Analytic Varieties, Advanced
  Studies in Pure Mathematics, vol.~1, 1983, pp.~219--230.

\end{thebibliography}

\newcommand{\etalchar}[1]{$^{#1}$}
\def\cprime{$'$} \def\cprime{$'$} \def\cprime{$'$} \def\cprime{$'$}
  \def\cprime{$'$} \def\dbar{\leavevmode\hbox to 0pt{\hskip.2ex
  \accent"16\hss}d} \def\cprime{$'$} \def\cprime{$'$}
  \def\polhk#1{\setbox0=\hbox{#1}{\ooalign{\hidewidth
  \lower1.5ex\hbox{`}\hidewidth\crcr\unhbox0}}} \def\cprime{$'$}
  \def\cprime{$'$} \def\cprime{$'$} \def\cprime{$'$}
  \def\polhk#1{\setbox0=\hbox{#1}{\ooalign{\hidewidth
  \lower1.5ex\hbox{`}\hidewidth\crcr\unhbox0}}} \def\cdprime{$''$}
  \def\cprime{$'$} \def\cprime{$'$} \def\cprime{$'$} \def\cprime{$'$}
\providecommand{\bysame}{\leavevmode\hbox to3em{\hrulefill}\thinspace}
\providecommand{\MR}{\relax\ifhmode\unskip\space\fi MR }
\providecommand{\MRhref}[2]{%
  \href{http://www.ams.org/mathscinet-getitem?mr=#1}{#2}
}
\providecommand{\href}[2]{#2}

\bigskip

  Princeton University, Princeton NJ 08544-1000, \

\email{kollar@math.princeton.edu}

\end{document}